\documentclass{article}
\usepackage{amsthm,amsfonts,amsmath,amssymb}

\newcommand{\eq}{\begin{equation}}
\newcommand{\en}{\end{equation}}

\newcommand{\prob}{\mathbb P}

\newcommand{\e}{{\rm e}}
\newcommand{\dd}{{\rm d}}
\newcommand{\htau}{\hat{\tau}}

\newtheorem{theorem}{\large Theorem}
\newtheorem{proposition}[theorem] {\large Proposition}

\begin{document}
\title{Winning Rate in the Full-Information Best Choice Problem}
\author{Alexander V. Gnedin\thanks{Postal address:
 Department of Mathematics, Utrecht University,
 Postbus 80010, 3508 TA Utrecht, The Netherlands. E-mail address: gnedin@math.uu.nl}~~~
and~~~ 
Denis I. Miretskiy \thanks{Postal address:
 Department of Applied Mathematics, University of Twente,
 Postbus 217, 7500 AE Enschede, The Netherlands. E-mail address: d.miretskiy@math.utwente.nl}}
\date{}
\maketitle

\begin{abstract}
\noindent
Following a long-standing suggestion by Gilbert and Mosteller,
we derive an explicit formula for the asymptotic winning rate in the full-information problem of the best choice.
\end{abstract}

\section{Introduction}

Let $X_1,X_2\ldots$  be a sequence of independent  uniform $[0,1]$ random variables.
The full-information best choice  problem, as introduced
by Gilbert and Mosteller \cite{GM}, 
asks one to find a stopping rule $\tau_n$  to maximise the probability
\eq\label{stop}
P_n(\tau):=\prob(X_{\tau}=\max(X_1,\ldots,X_n))
\en
over all stopping rules $\tau\leq n$ adapted to the sequence $(X_i)$.
The name `full information' was attached to the problem to stress that the observer
learns the exact values of $X_i$'s and knows their distribution, in contrast 
to the `no information' problem where only the relative ranks of observations are available
(see \cite{SamSurv} for a survey and history of the best choice or `secretary' problems).
Because  the stopping criterion (\ref{stop}) depends only on ranks of  the observations,
the instance of uniform distribution covers, in fact, the  general case of sampling from arbitrary
continuous distribution.

\par Gilbert and Mosteller showed that
the optimal stopping rule is
of the form
$$\tau_n=\min\{i:\,X_i=\max(X_1,\ldots,X_i)~{\rm and~}X_i\geq d_{n-i}\},$$
where  $d_k$ is a sequence of  decision numbers defined by the equation
\eq\label{d}
\sum_{j=1}^k (d_k^{-j}-1)/j=1 {~~\rm for~~}k\geq 1,~~{\rm and~~}d_0=0.
\en
They also  proved that $d_k\uparrow 1$ in such a way that
$k(1-d_k)\to c$ for $c=0.804\ldots$ the solution to the transcendental equation 
\eq\label{c}
\int_0^c x^{-1}({e^x-1})\,{\rm d}x=1\,,
\en
and they provided a numerical evidence that the optimal probability of the best choice
$P_n^*:=P_n(\tau_n)$
converges to a limit $P^*=0.580164\ldots$ The limiting value was justified by different methods 
in the subsequent work \cite{BG, GnFI, SamUn, SamSurv}
 along with the explicit formula 
\eq\label{limpr}
P^*=\e^{-c}+(\e^c-c-1)\int_1^\infty \e^{-cx}x^{-1}\,\dd x\,
\en
due to Samuels \cite{SamUn}.

\par Refinements and generalisations of the results of \cite{GM} appeared in 
 \cite{GnFI, GnBC, payoff, HillKennedy, Pet, Por1, Sam-Cahn}. Still,
one interesting feature of the optimal stopping rule  seems to have not been 
discussed in the literature. 
We mean the tiny
Section 3e in \cite{GM} where Gilbert and Mosteller say: 
``One would correctly anticipate that as $n$ increases, the probability of winning at a given draw tends to zero.
On the other hand,
$\,n\,\prob({\rm win~at~draw}~i)$ tends to a constant for 
$i/n$ tending to a constant $\lambda\,$". 
 Spelled out in detail, Gilbert and Mosteller claimed  existence of the limit
\eq\label{conv}
w(t)=\lim_{i,n\to\infty,\,i/n\to t}
 n\,\prob(\tau_n=i,\,X_i=\max(X_1,\ldots,X_n))
\en
where $t\in [0,1]$ stands for their $\lambda$. Such a function may be called the asymptotic
{\it winning rate} since it tells us how the chance of correct recognising the maximum 
is distrubuted over the time, hence
the total probability of the best choice must satisfy 
$$P^*=\int_0^1 w(t)\,\dd t\,.$$
In this paper, we prove the conjecture of \cite{GM} regarding the convergence 
and we derive an explicit formula for the winning rate (\ref{conv}). In fact, we show more: the function $w$ appears as the 
exact winning rate in a continuous-time
version of the best choice problem associated with a planar Poisson process, as developed in 
\cite{GnFI, GnBC, Sam3}.

\section{The Poisson framework}
We start by recalling the setup from \cite{GnFI, GnBC}. Consider a homogeneous planar Poisson process (PPP) in the
semi-infinite strip
$R=[0,1]\times \,]-\infty,0]$, with Lebesgue measure as intensity. 
The generic atom $a=(t,x)\in R$ of the PPP is understood as score $x$
observed at time $t$ . 
Let ${\cal F}=({\cal F}_t, ~t\in [0,1])$ be the filtration with ${\cal F}_t$ the $\sigma$-algebra generated by 
the PPP restricted to $[0,t]\times  \,]-\infty,0]$. 
We say that an atom $a=(t,x)$ of the PPP is a {\it record} if there are no 
other PPP-atoms north-west of $a$. The maximum of the PPP is an atom $a^*=(t^*,x^*)$ with the largest $x$-value. 
Alternatively, the maximum $a^*$ can  be defined as the last record of the PPP, that is the record with the largest $t$-value.
For $\tau$ a ${\cal F}$-adapted stopping rule with values in $[0,1]$, 
the performance of $\tau^*$ is defined as the probability of 
the event $\{\tau=t^*\}$, interpreted as the best choice from the PPP. 
The associated best choice problem amounts to maximising probability of the event $\{\tau=t^*\}$.

\par In a poissonised version of the Gilbert-Mosteller problem the observations sampled from 
the $[0,1]$ uniform distribution arrive 
on $[0,\ell]$ at epochs of a rate 1 Poisson process \cite{BG, Bojd, GnSak, Sak}. This is equivalent to 
the PPP setup with background space $[0,\ell]\times[0,1]$, which can be 
mapped  linearly onto $[0,1]\times[-\ell,0]$ so that the componentwise order of points is preserved.
Now, the optimal stopping in 
$[0,1]\times[-\ell,0]$
fits in the framework with the background space $R$
by a minor modification of the  
stopping criterion: 
a  stopping rule $\tau$  adapted to $\cal F$ is evaluated 
by the probability of the event $\{\tau=t^*, a^*>-\ell\}$ that stopping occurs at the maximum atom 
and above $-\ell$. In this sense we shall speak of a {\it constrained} best choice problem.

\par Let $\Gamma=\{(t,x)\in R: -x(1-t)<c\}$ where $c$ is as in (\ref{c}). 
It is known \cite{GnFI} that the optimal stopping rule is the first time (if any)
when the record process enters $\Gamma$, that is
$$\tau^*=\min\{t: {\rm~ there~is~a~record~}a=(t,x)\in \Gamma\}$$
(or $\tau^*=1$ if no such $t\in [0,1[$ exists). 
Similarly,  the optimal stopping rule for the constrained problem is the first time (if any)
when the record process enters $\Gamma(\ell):=\Gamma\cap ([0,1]\times[0,-\ell])$.
\par Let
$$g(\ell,t):=\prob(\tau^*=t^*,\, t^*<t,\, x^*>-\ell)$$
be the probability that $\tau^*$ wins by stopping above $-\ell$ and before $t$ and let 
$$g(\infty,t):=\prob(\tau^*=t^*,\, t^*<t).$$
By the above relation between the constrained and unconstrained problems we have
$$g(\ell,t)=g(\infty,t)~~~{\rm for~~~}0\leq t\leq (1-c/\ell)_+\,.$$  
The {\it winning rate} in the Poisson problem is defined as
$$w(t)=\partial_t\, g(\infty,t).$$

\section{Computing the rate}

Because the atoms south-west of $(x,-\ell)$ fall outside the stopping region $\Gamma(\ell)$ we have
$\partial_\ell\,g(\ell,t)=0$ and
$g(\infty,t)=g(c/(1-t),t)$ for $\ell>s/(1-t)$. To determine $\partial_\ell\,g(\ell,t)$ for $\ell>c/(1-t)$
consider two rectangles $R_1=[0,t]\times[-\ell,0]$ and $R_2=[0,t]\times[-\ell+\delta,0]$ with small $\delta>0$.
The optimal constrained stopping rules in $R_1$ and $R_2$ stop before $t$ at distinct atoms if and only if
the record process enters $\Gamma({\ell})$ at some atom $a_0=(\sigma,\xi)\in [(1-c/\ell)_+\,,t]\times   
[-\ell,-\ell+\delta]$.  Then stopping at $a_0\in R_1\setminus R_2$ is a win if $a_0=a^*$, which occurs with probability
$$p_1={\rm e}^{-\ell}\ell(t-(1-c/\ell)_+){\delta\over\ell}+o(\delta)=
{\rm e}^{-\ell}(t-(1-c/\ell)_+){\delta}+o(\delta).$$
On the other hand, stopping in $R_2$ is a win (and stopping at $a_0$ is a loss) if $a_0$ is folowed by some
$k>0$ atoms in $[\sigma,1]\times[-\ell,0]$, the leftmost of these $k$ atoms appears within $[\sigma,t]\times[-\ell,0]$
and it is the overall maximum $a^*$ 
which is an event of probability
$$p_2={\rm e}^{-\ell}\sum_{k=1}^\infty{c^{k+1}\over(k+1)!}\left[
1-(k+1){(t-(1-c/\ell)_+)\over c/\ell}\,\,{(1-t)^k\over (c/\ell)^k}-{(1-t)^{k+1}\over (c/\ell)^{k+1}}
\right]{1\over k}\,{\delta\over \ell}+o(\delta).$$
It follows that 
$$\partial_\ell\,g(\ell,t)=\lim_{\delta\to 0}{p_1-p_2\over \delta}={\rm e}^{-\ell}\int\limits_{(1-c/\ell)_+}^t
\left(
1-\sum_{k=1}^\infty\left[{\ell^k(1-\sigma)^k\over k!\,k}-{\ell^k(1-t)^k\over k!\,k}
\right]
\right){\rm d}\sigma\,.$$   
Now, computing the mixed second derivative $\partial_{\ell\,t} \,g(\ell,t)$ and integrating 
in $\ell$ from $0$  to $c/(1-t)$ we obtain the winning rate in the Poisson problem, which is our main result.

\begin{proposition} The winning rate is given by the formula
\eq\label{w}
w(t)=-{\rm e}^{-c}+{\e^{-ct}-\e^{-ct/(1-t)}\over t}+
{e^{-c t}-te^{-c}\over 1-t}+
{c\over 1-t}
\left\{ I\left({c\over1-t}\,,\,c\right)-I\left({ct\over 1-t}\,,\,ct\right)\right\} ,
 \en
where $c$ is as in {\rm (\ref{c})} and for $0<s<t$
$$I(t,s)=\int_s^t \xi^{-1}{\e^{-\xi}}\,{\rm d}\xi\,.$$
 \end{proposition}
\noindent
The boundary values of $w$ are $w(0)=1-\e^{-c}=0.5526\ldots$ and $w(1)=\e^{-c}=0.4473\ldots$,
in accordance with \cite[Fig. 3]{GM}.
A {\tt Mathematica}-drawn graph of {\rm (\ref{w})} exhibits a  curve identical to that in 
\cite[Fig 3.]{GM}.

\par The special value (\ref{c}) of $c$ was not used in the argument, hence 
the right side of the formula (\ref{w}) gives the winning rate for 
every stopping rule defined 
by a stopping region  like $\Gamma$ but with arbitrary positive value of the constant in place of $c$.
We also note that
the winning rate in the constrained problem coincides with $w(t)$ for $t<(1-c/\ell)_+$.


\section{Embedding and convergence}

It remains to show that  $w$ given by (\ref{w}) is indeed the limiting value for the finite-$n$ problem 
in (\ref{conv}).
To that end, we will exploit the embedding  technique from \cite{GnFI}.

\par With $n$ fixed, divide $R$ in strips $J_i=[(i-1)/n, \,i /n[\,\,\times\,\,]\!-\infty,0]$,
$i=1,\ldots,n$. Consider a sequence 
${\cal Y}_n=((T_i,Y_i), \,i=1,\ldots,n)$ where $(T_i,Y_i)$ is an atom with the largest 
$x$-component within the strip $J_i$. 
Observe that the point process of records in ${\cal Y}_n$ is 
a subset of the set of records of the PPP in $R$,
in particular $\max Y_i=x^*$.
By  homogeneity of the PPP we have
$T_1,Y_1,\ldots,T_n,Y_n$ 
jointly independent, with each $T_i$  uniformly distributed on $[(i-1)/n, \,i /n[$
and each $Y_i$ exponentially distributed on $]-\infty,0]$ with rate $1/n$. 
It follows that the discrete-time optimal stopping problem of recognising the maximum in ${\cal Y}_n$ 
is equivalent to the Gilbert-Mosteller problem with exponentially distributed observations.

\par Let $\htau_n$ be the optimal stopping rule for recognising the maximum in ${\cal Y}_n$.
We shall view $\htau_n$ as a strategy 
for choosing the maximum of PPP with the additional option of 
{\it partial return} meaning that $\htau_n$ assumes values in $[0,1]$, that
$$\{(i-1)/n< \htau_n\leq i/n\}\in {\cal F}_{i/n}$$
and that $\{(i-1)/n< \htau_n\leq i/n\}$ is associated with the stopping at $(T_i,Y_i)$.
Explicitly, $\htau_n$ stops at the first time the sequence of ${\cal Y}_n$-records enters
$$\Gamma_n=\bigcup\limits_{i=1}^n ~](i-1)/n,i/n]\times [b_{n-i},0],$$
where $b_k=n\log d_k$ and the $d_k$'s are the decision numbers as in  (\ref{d}).
The partial return option 
implies that the winning chance of $\htau_n$ is higher than that of
$\tau^*$.

\par Let $a'$ be the last record before $a^*$. One checks easily that $\tau^*$ and $\htau$ may differ
only if either $a^*$ or $a'$ hit the domain
$$\Delta_n:=(\Gamma_n\setminus\Gamma) \cup(\Gamma\setminus\Gamma_n).$$ 
By \cite[Equation (11)]{GnFI} we have $((i-1)/n,b_{n-i})\not\in\Gamma$ and
$(i/n,b_{n-i})\in\Gamma$ for $i=1,\ldots,n$. This combined with the fact that 
the distribution of $t^*$ is uniform and that of $x^*$ is exponential yields
$$n\,\prob(a^*\in \Delta_n\cap J_i)<\exp\left({-nc\over n-i+1}\right)-\exp\left({-nc\over n-i}\right)=O(n^{-1})$$
uniformly in $i\leq n$. A similar estimate holds also for $a'$, and because

$$\prob((i-1)/n<\htau_n\leq i/n)=\prob(\tau_n=i), ~~~w(i/n)=\prob(\tau^*=t^*, a^*\in J_i)+O(n^{-1})$$
(the second since $w$ is smooth on $[0,1]$)
we conclude:

\begin{proposition} As $n\to\infty$ the optimal stoping rule $\tau_n$ satisfies
$$
\max_{1\leq i\leq n}|w(i/n)-n\prob(\tau_n=i,~X_{i}=\max(X_1,\ldots,X_n)|=O(n^{-1}).
$$
where $w$ is given by {\rm (\ref{w})}.
\end{proposition}

\end{document}